\documentclass[12pt,leqno]{amsart}
\usepackage{amssymb,verbatim,enumerate,ifthen,cite}
\usepackage[mathscr]{eucal}
\usepackage[utf8]{inputenc}
\usepackage[T1]{fontenc}
\def\N{\mathbb{N}}
\def\R{\mathbb{R}}

\def\F{\mathscr{F}}
\def\C{\mathscr{C}}

\def\M{\mathscr{M}}
\def\A{\mathscr{A}}
\def\H{\mathscr{H}}
\def\E{{\mathscr{E}}}

\newtheorem{theorem}{Theorem}[section]
\newtheorem*{theorem*}{Theorem}
\def\Thm#1#2{\ifthenelse{\equal{#1}{*}}{\begin{theorem*}#2\end{theorem*}}
             {\begin{theorem}\label{T#1}#2\end{theorem}}}
\newtheorem{Atheorem}{Theorem}

\def\thm#1{Theorem~\ref{T#1}}
\newtheorem{proposition}[theorem]{Proposition}
\newtheorem*{proposition*}{Proposition}
\def\Prp#1#2{\ifthenelse{\equal{#1}{*}}{\begin{proposition*}#2\end{proposition*}}
{\begin{proposition}\label{P#1}#2\end{proposition}}}
\def\prp#1{Proposition~\ref{P#1}}

\newtheorem{corollary}[theorem]{Corollary}
\newtheorem*{corollary*}{Corollary}
\def\Cor#1#2{\ifthenelse{\equal{#1}{*}}{\begin{corollary*}#2\end{corollary*}}
             {\begin{corollary}\label{C#1}#2\end{corollary}}}
\def\cor#1{Corollary~\ref{C#1}}

\newtheorem{lemma}[theorem]{Lemma}
\newtheorem*{lemma*}{Lemma}
\def\Lem#1#2{\ifthenelse{\equal{#1}{*}}{\begin{lemma*}#2\end{lemma*}}
             {\begin{lemma}\label{L#1}#2\end{lemma}}}
\def\lem#1{Lemma~\ref{L#1}}

\theoremstyle{definition}
\newtheorem{remark}[theorem]{Remark}
\newtheorem*{remark*}{Remark}
\def\Rem#1#2{\ifthenelse{\equal{#1}{*}}{\begin{remark}\rm #2\end{remark}}
             {\begin{remark}\label{R#1}\rm #2\end{remark}}}

\newtheorem{example}[theorem]{Example}
\newtheorem*{example*}{Example}
\def\Exa#1#2{\ifthenelse{\equal{#1}{*}}{\begin{example*}\rm #2\end{example*}}
             {\begin{example}\label{Ex#1}\rm #2\end{example}}}

\def\eq#1{{\rm(\ref{E#1})}}
\def\Eq#1#2{\ifthenelse{\equal{#1}{*}}
  {\begin{equation*}\begin{aligned}#2\end{aligned}\end{equation*}}
  {\begin{equation}\begin{aligned}\label{E#1}#2\end{aligned}\end{equation}}}

\begin{document}
\begin{flushright}
\end{flushright}
\vspace{5mm}

\date{\today}

\title{On approximately convex and affine functions}

\author[A. R. Goswami]{Angshuman R. Goswami}
\address[A. R. Goswami]{Doctoral School of Mathematical and Computational Sciences, University of Debrecen, H-4002 Debrecen, Pf.\ 400, Hungary}
\email{angshu@science.unideb.hu}
\author[Zs. P\'ales]{Zsolt P\'ales}
\address[Zs. P\'ales]{Institute of Mathematics, University of Debrecen, 
H-4002 Debrecen, Pf.\ 400, Hungary}
\email{pales@science.unideb.hu}


\thanks{The research of the second author was supported by the K-134191 NKFIH Grant and the 2019-2.1.11-TÉT-2019-00049, the EFOP-3.6.1-16-2016-00022 and the EFOP-3.6.2-16-2017-00015 projects. The last two projects are co-financed by the European Union and the European Social Fund.}

\begin{abstract}
A real valued function $f$ defined on a real open interval $I$ is called $\Phi$-convex if, for all $x,y\in I$, $t\in[0,1]$ it satisfies
\Eq{*}{
  f(tx+(1-t)y)\leq tf(x)+(1-t)f(y)+t\Phi\big((1-t)|x-y|\big)+(1-t)\Phi\big(t|x-y|\big),
}
where $\Phi:\R_+\to\R_+$ is a  nonnegative error function. If $f$ and $-f$ are simultaneously $\Phi$-convex, then $f$ is said to be a $\Phi$-affine function. In the main results of the paper, we describe the structural and inclusion properties of these two classes. We characterize these two classes of functions and investigate their relationship with approximately monotone and approximately-H\"older functions. We also introduce a subclass of error functions which enjoy the so-called $\Gamma$ property and we show that the error function which is the most optimal for a $\Phi$-convex function has to belong to this subclass. The properties of this subclass of error function are investigated as well. Then we offer two formulas for the lower $\Phi$-convex envelop. Besides, a sandwich type theorem is also added. 
\end{abstract}

\maketitle

\section{Introduction} 

The motivation for our paper comes from the theory of approximate convexity which has become an active field of research and many important contributions have been made, see for instance \cite{BorNag13,BurHaz11,BurHazJuh11,DilHowRob99,DilHowRob00,DilHowRob02,Ger88c,GilGonNikPal17,GonNikPalRoa14,Gre52,HyeUla52,Haz05a,Haz10,HazPal04,HazPal05,JarLac09,JarLac10,Krz01,Lac99,Mak17,MakHaz17,MakPal10b,MakPal11a,MakPal12b,MakPal12c,MakPal12a,MakPal13a,MakPal13b,Mro01,MroTabTab09,NgNik93,NikPal01,NikPal07,Pal01c,Pal03a,Rol79,SpuTab12,TabTabZol10b,TabTabZol10a}.
In these papers several aspects of approximate convexity were investigated: stability problems, Bernstein--Doetsch-type theorems, Hermite--Hadamard type inequalities, etc. The notions and results have various applications in nonsmooth and convex analysis and optimization theory, and also in the theory of functional equations and inequalities. The main concepts and results of this paper are distillated from the following elementary observations; primarily studied in the paper \cite{MakPal12a}. 

Assume that $I$ is a nonempty interval and a function $f:I\to\R$ satisfies the following inequality 
\Eq{0}{
  f(tx+(1-t)y)\leq tf(x)+(1-t)f(y)+\varepsilon\big(t\big((1-t)|x-y|\big)^{p}+(1-t)\big(t|x-y|\big)^{p}\big), \\
\qquad (x,y\in I,\,t\in[0,1])
}
for some nonnegative constant $\varepsilon$ and real constant $p\in\R$. In other words, $f$ is approximately convex with an error term described in terms of the $p$th power function. Clearly, if $\varepsilon=0$, then the above condition is equivalent to the convexity of $f$. One can easily notice that every convex function $f$ satisfies \eq{0}. On the other hand, if $p=0$, then $f$ satisfies \eq{0} if and only if it can be decomposed as $f=g+h$, where $g$ is convex and $||h||_\infty\leq\varepsilon/2$. In the case $p=1$ the results of the paper \cite{Pal03a} show that $f$ fulfills \eq{0} if and only if $f=g+\ell$, where $g$ is convex and $\ell$ is Lipschitzian with a Lipschitz modulus $c\varepsilon$. Surprisingly, for $p>2$, the situation is completely different. Then \eq{0} holds if and only if $f$ is convex.

The above described observations and results motivate the investigation of classes of functions that obey a more general approximate convexity and also the related approximately affine property. In fact, the class of $\Phi$-convex functions was introduced in the paper \cite{MakPal12a}, but this property was only characterized therein. In this work, we describe the structural properties of this function class and determine the error function which is the most optimal one. We show that optimal error functions for approximate convexity must possess the so-called $\Gamma$ property. Then we offer a precise formula for the $\Phi$-convex envelope and also obtain sandwich-type theorems.

\section{Basic Results}

Let $I$ be a nonempty open real interval throughout this paper and let $\ell(I)\in\,]0,\infty]$ denote its length. The symbols $\R$ and $\R_+$ denote the sets of real and nonnegative real numbers, respectively. 

The class of all functions $\Phi:[0,\ell(I)[\,\to\R_+$, called error functions, will be denoted by $\E(I)$. Obviously, $\E(I)$ is a convex cone, i.e., it is closed with respect to addition and multiplication by nonnegative scalars. In what follows, we are going to define four properties related to an error function $\Phi\in\E(I)$. First we recall the notions of $\Phi$-monotone and $\Phi$-H\"older functions that have been introduced in our former papers \cite{GosPal20,GosPal21}.

A function $f:I\to\R$ will be called \emph{$\Phi$-monotone} if, for all $x,y\in I$ with $x\leq y$, 
\Eq{H1}{
  f(x)\leq f(y)+\Phi(y-x).
}
If this inequality is satisfied with the identically zero error function $\Phi$, then we say that $f$ is \emph{monotone (increasing)}.
The class of $\Phi$-monotone functions on $I$ will be denoted by $\M_\Phi(I)$.

A function $f:I\to\R$ will be called \emph{$\Phi$-H\"older} if, for all $x,y\in I$, 
\Eq{H2}{
  |f(x)-f(y)|\leq\Phi(|x-y|).
}
The class of $\Phi$-Hölder functions on $I$ will be denoted by $\H_\Phi(I)$.

The following concept was introduced in the paper \cite{MakPal12a}. We say that a function $f:I\to\R$ is \emph{$\Phi$-convex} if, for all $x,y\in I$ and $t\in[0,1]$, 
\Eq{1}{
  f(tx+(1-t)y)\leq tf(x)+(1-t)f(y)+t\Phi\big((1-t)|x-y|\big)+(1-t)\Phi\big(t|x-y|\big).
}
If the above inequality is satisfied with the identically zero error function $\Phi$, then we say that $f$ is \emph{convex}.
The class of $\Phi$-convex functions on $I$ will be denoted by $\C_\Phi(I)$.

Finally, a function $f:I\to\R$ is said to be \emph{$\Phi$-affine} if, for all $x,y\in I$ and $t\in[0,1]$, 
\Eq{2}{
  |f(tx+(1-t)y)-tf(x)-(1-t)f(y)|
  \leq t\Phi\big((1-t)|x-y|\big)+(1-t)\Phi\big(t|x-y|\big).
}
If, in particular, the above inequality is satisfied with the identically zero error function $\Phi$, then we say that $f$ is \emph{affine}.
The class of $\Phi$-affine functions on $I$ will be denoted by $\A_\Phi(I)$.

\Prp{MC1}{Let $\Phi_1,\dots,\Phi_n\in\E(I)$ and $\alpha_1,\dots,\alpha_n\in\R_+$. Then
\Eq{*}{
  \alpha_1\C_{\Phi_1}(I)+\cdots+\alpha_n\C_{\Phi_n}(I)
  \subseteq \C_{\alpha_1\Phi_1+\cdots+\alpha_n\Phi_n}(I).
}
In particular, for all functions $\Phi\in\E(I)$, the class $\C_\Phi(I)$ is convex.}

\begin{proof}
To prove the inclusion of the proposition, let $f\in \alpha_1\C_{\Phi_1}(I)+\cdots+\alpha_n\C_{\Phi_n}(I)$. Then, $f$ is of the form 
\Eq{999}{
  f=\alpha_1 f_1+\dots+\alpha_n f_n
}
for some elements $f_1,\dots,f_n$ belonging to $\C_{\Phi_1}(I),\dots,\C_{\Phi_n}(I)$, respectively. Then, for all $x,y\in I$, $t\in[0,1]$, and $i\in\{1,\dots,n\}$, we have
\Eq{*}{
  f_i(tx+(1-t)y)\leq tf_i(x)+(1-t)f_i(y)+t\Phi_i\big((1-t)|x-y|\big)+(1-t)\Phi_i\big((t|x-y|\big).
}
Multiplying this inequality by $\alpha_i$ and summing up side by side, we will arrive at the inequality \eq{1}, where $\Phi:=\sum_{i=1}^{n}\alpha_i\Phi_i$.
This proves that $f\in\C_{\Phi}(I)$ and completes the proof of the second inclusion.
\end{proof}

The following result is the counterpart of the previous statement.

\Prp{HA1}{Let $\Phi_1,\dots,\Phi_n\in\E(I)$ and $\alpha_1,\dots,\alpha_n\in\R$. Then
\Eq{*}{
  \alpha_1\A_{\Phi_1}(I)+\cdots+\alpha_n\A_{\Phi_n}(I)
&\subseteq \A_{|\alpha_1|\Phi_1+\cdots+|\alpha_n|\Phi_n}(I).
}
In particular, for all functions $\Phi\in\E(I)$, the class $\A_\Phi(I)$ is convex and central symmetric, i.e., $\A_\Phi(I)$ is closed with respect to multiplication by $(-1)$.}

\begin{proof}
To prove the inclusion of the proposition, let $f\in \alpha_1\A_{\Phi_1}(I)+\cdots+\alpha_n\A_{\Phi_n}(I)$. Then, $f$ is of the form \eq{999} for some elements $f_1,\dots,f_n$ belonging to $\A_{\Phi_1}(I),\dots,\A_{\Phi_n}(I)$, respectively. Then, for all $x,y\in I$, $t\in[0,1]$, and $i\in\{1,\dots,n\}$, we have
\Eq{*}{
  |f_i(tx+(1-t)y)-tf_i(x)-(1-t)f_i(y)|
  \leq t\Phi_i\big((1-t)|x-y|\big)+(1-t)\Phi_i\big((t|x-y|\big).
}
Multiplying this inequality by $|\alpha_i|$ and summing up side by side, we will arrive at the inequality 
\Eq{*}{
  |f(tx+(1-t)y)-tf(x)-(1-t)f(y)|
  &=\bigg|\sum_{i=1}^n \alpha_i\big(f_i(tx+(1-t)y)-tf_i(x)-(1-t)f_i(y)\big)\bigg|\\
  &\leq\sum_{i=1}^n |\alpha_i|\big|f_i(tx+(1-t)y)-tf_i(x)-(1-t)f_i(y)\big|\\
  &\leq\sum_{i=1}^n |\alpha_i|t\Phi_i\big((1-t)|x-y|\big)+(1-t)\Phi_i\big((t|x-y|\big)\\
  &=t\Phi\big((1-t)|x-y|\big)+(1-t)\Phi\big((t|x-y|\big),
}
where $\Phi:=\sum_{i=1}^{n}|\alpha_i|\Phi_i$.
This proves that \eq{2} holds, i.e., $f\in\A_{\Phi}(I)$ and completes the proof of the second inclusion.
\end{proof}

\Prp{H0}{Let $\Phi\in\E(I)$. Then
\Eq{H0}{
  \H_\Phi(I)\subseteq\A_\Phi(I) \qquad\mbox{and}\qquad
\A_\Phi(I)=\C_\Phi(I)\cap(-\C_\Phi(I)).
}
}

\begin{proof}
To prove the inclusion, let $f\in\H_\Phi(I)$. Let $x,y\in I$ and $t\in[0,1]$. 
By the $\Phi$-H\"older property of $f$, we have 
\Eq{*}{
  f(tx+(1-t)y)&\leq f(x)+\Phi\big((1-t)|x-y|\big), \\
  f(tx+(1-t)y)&\leq f(y)+\Phi\big(t|x-y|\big).
}
Multiplying the above inequalities by $t$ and $1-t$, then summing up the inequalities so obtained side by side, we get
\Eq{*}{
  f(tx+(1-t)y)
  \leq tf(x)+(1-t)f(y)+t\Phi\big((1-t)|x-y|\big)+(1-t)\Phi\big((t|x-y|\big),   
}
which shows that $f\in\C_\Phi(I)$. Repeating the same argument with $(-f)$ instead of $f$, it follows that $-f\in\C_\Phi(I)$. Therefore, $f\in\C_\Phi(I)\cap(-\C_\Phi(I))=\A_\Phi(I)$.

To show the equality assertion, assume that $f$ is a $\Phi$-affine function. Then, for any $x,y\in I$ and $t\in[0,1]$, $f$  will satisfy the inequality \eq{2} and hence the inequalities
\Eq{Mon}{
  f(tx+(1-t)y)\leq &tf(x)+(1-t)f(y)+t\Phi\big((1-t)|x-y|\big)+(1-t)\Phi\big(t|x-y|\big).\\
\mbox{and}\\
tf(x)+(1-t)f(y)\leq& f(tx+(1-t)y)+t\Phi\big((1-t)|x-y|\big)+(1-t)\Phi\big(t|x-y|\big).
}
holds. Rearranging theses two inequalities, we have that both $f$ and $-f$ are $\Phi$-affine. That is
 $f\in \C_\Phi(I)\cap(-\C_\Phi(I)).$

To show the inverse inclusion, let $f\in \C_\Phi(I)\cap(-\C_\Phi(I)).$  Due to the property of $\Phi$-convexity of the two classes of function, $f$ will satisfy the two inequalities in \eq{Mon}. Hence, inequality \eq{2} holds for any $x, y\in I$ and $t\in[0,1]$. This inequality being symmetric in $x$ and $y$, we get that \eq{2} is satisfied for all $x,y\in I$.
\end{proof}

\Cor{0}{Let $\Phi,\Psi\in\E(I)$. Then
\Eq{*}{
  \H_\Phi(I)+\C_\Psi(I)\subseteq \C_{\Phi+\Psi}(I)
  \qquad\mbox{and}\qquad
  \H_\Phi(I)+\A_\Psi(I)\subseteq \A_{\Phi+\Psi}(I)
}}

\begin{proof} By \prp{MC2}, we have $\H_\Phi(I)\subseteq\A_\Phi(I)\subseteq\C_\Phi(I)$. Therefore, the statement follows from \prp{MC1}.
\end{proof}

We say that a family $\F$ of real valued functions is \emph{closed with respect to the pointwise supremum} if $\{f_\gamma:I\to\R\mid\gamma\in\Gamma\}$ is a subfamily of $\F$ with a pointwise supremum $f:I\to\R$, i.e., 
\Eq{sup}{
  f(x)=\sup_{\gamma\in\Gamma} f_\gamma(x)\qquad(x\in I),
}
then $f\in\F$. Similarly, we can define that a family $\F$ of real valued functions is \emph{closed with respect to the pointwise infimum}. A family $\{f_\gamma:I\to\R\mid\gamma\in\Gamma\}$ is called a \emph{chain} if, for all $\alpha,\beta\in\Gamma$, either $f_\alpha\leq f_\beta$ or $f_\beta\leq f_\alpha$ holds on $I$. We say that a family $\F$ of real valued functions is \emph{closed with respect to the pointwise chain supremum (chain infimum)} if $\{f_\gamma:I\to\R\mid\gamma\in\Gamma\}\subseteq\F$ is a chain with a pointwise supremum (infimum) $f:I\to\R$, then $f\in\F$. 

The following result was established in the paper \cite{MakPal12a}, therefore its proof is omitted.

\Prp{MC2}{Let $\Phi\in\E(I)$. Then the class $\C_\Phi(I)$ is closed under pointwise supremmum and  pointwise chain infimum.}

\Cor{MC2}{Let $\Phi\in\E(I)$. Then the class $\C_\Phi(I)$ is closed with respect to the limsup operation.}

\begin{proof}
Assume that $f:I\to\R$ is the upper limit of a sequence $f_n:I\to\R$. Then
\Eq{*}{
  f=\inf_{n\in\N} g_n, \qquad\mbox{where}\qquad g_n:=\sup_{k\geq n} f_k.
}
If all the functions $f_n$ are $\Phi$-convex, then 
for all $n\in\N$, the function $g_n$ is $\Phi$-convex.
(Because these classes are closed with respect to the pointwise supremum.) On the other hand, the sequence $(g_n)$ is decreasing, therefore $f$ is the pointwise chain infimum of $\{g_n\mid n\in\N\}$, thus $f$ is also $\Phi$-convex. 
\end{proof}

\Prp{H2}{Let $\Phi\in\E(I)$. Then the class $\A_\Phi(I)$ is closed under pointwise chain infimum and pointwise chain supremum. Consequently $\A_\Phi(I)$ is closed with respect to the pointwise limit operation.}

\begin{proof} 
The statement easily follows from \prp{MC2}. The statement concerning the pointwise limit operation is obvious.
\end{proof}

\Thm{Com}{Let $\Phi\in\E_0(I)$, let $f:I \to \R$ be a $\Phi$-convex function and $g:\R \to \R_+ $ be an increasing and sublinear function. Then $(g\circ f)$ is a $(g\circ\Phi)$-convex function.}

\begin{proof}
Using the $\Phi$-convexity of $f$, the nondecreasingness and sublinearity and properties of $g$, we arrive at
\Eq{*}{
(g\circ f)(tx+(1-t)y)
&\leq g\Big(t \big(f(x)+\Phi\big((1-t)|y-x|\big)\big)+(1-t)\big(f(y)+\Phi\big(t|y-x|\big)\big)\Big)\\
&\leq t \big(g\circ f(x)+g\circ \Phi\big((1-t)|y-x|\big)\big)+(1-t)\big(g\circ f(y)+g\circ \Phi\big(t|y-x|\big)\big).
}
Therefore, we arrive at our desired conclusion.
\end{proof}

\section{Characterizations of $\Phi$-convex and $\Phi$-affine functions}

The equivalence of assertions (i), (ii) and (iii) of the following result was established in the paper \cite{MakPal12a}. For the convenience of the reader, we provide the complete proof.

\Thm{1}{Let $\Phi\in\E(I)$ and $f:I\to\R$. Then the following conditions are equivalent to each other:
\begin{enumerate}[(i)]
 \item $f$ is $\Phi$-convex;
 \item For all $x,u,y\in I$ with $x<u<y$, 
 \Eq{*}{
 \dfrac{f(u)-f(x)-\Phi(u-x)}{u-x}\leq \dfrac{f(y)-f(u)+\Phi(y-u)}{y-u};
 }
 \item There exists a function $\varphi: I\to\R$ such that, for all $x,u\in I$,
 \Eq{ux}{
 f(u)+(x-u)\varphi(u)\leq f(x)+\Phi(|u-x|);
 }
 \item For all $n\in\N$, $x_1,\dots,x_n\in I,$ $t_1,\dots,t_n\geq0$ with $ t_1+\dots+t_n=1$,
 \Eq{*}{
 f(t_1x_1+\dots+t_nx_n)\leq \sum_{i=1}^{n}t_i\big(f(x_i)+\Phi(|(t_1x_1+\dots+t_nx_n)-x_i|)\big).
 }
\end{enumerate}
}

\begin{proof} {\it(i)$\Rightarrow$(ii):} Assume that $f$ is $\Phi$ convex and let $x<u<y$ be arbitrary elements of $I$. Choose $t\in[0,1]$ such that $u=tx+(1-t)y$. That is let $t=\dfrac{y-u}{y-x}$. Then \eq{1} can be rewritten as
\Eq{*}{
  f(u)\leq\dfrac{y-u}{y-x}f(x)+\dfrac{u-x}{y-x}f(y)
      +\dfrac{y-u}{y-x}\Phi(u-x)+\dfrac{u-x}{y-x}\Phi(y-u).
}
Therefore, 
\Eq{*}{
  (y-u+u-x)f(u)\leq&(y-u)f(x)+(u-x)f(y)\\
  &+(y-u)\Phi(u-x)+(u-x)\Phi(y-u).
}
Rearranging this inequality, it follows that
\Eq{*}{
 \dfrac{f(u)-f(x)-\Phi(u-x)}{u-x}&\leq \frac{f(y)-f(u)+\Phi(y-u)}{y-u}.
}

{\it(ii)$\Rightarrow$(iii):}
Assume that {\it(ii)} holds and define the function $\varphi$ on $I$ by
\Eq{*}{
  \varphi(u):=\underset{x<u}{\mathrm{\sup}}\bigg(\dfrac{f(u)-f(x)-\Phi(u-x)}{u-x}\bigg) \qquad(u\in I).
}
In view of condition {\it(ii)}, for all $x<u<y$ in $I$, we have
\Eq{6}{
\dfrac{f(u)-f(x)-\Phi(u-x)}{u-x}\leq\varphi(u)\leq \dfrac{f(y)-f(u)+\Phi(y-u)}{y-u}.
}
From the left hand side inequality in \eq{6}, we get
\Eq{7x}{
f(u)+(x-u)\varphi(u)\leq f(x)+\Phi(u-x) \qquad (x\in I,\,x<u).
}
Similarly, from the right hand side inequality in \eq{6} (replacing $y$ by $x$, it follows that
\Eq{7y}{
f(u)+(x-u)\varphi(u)\leq f(x)+\Phi(x-u) \qquad (x\in I,\,u<x).
}
Now, combining inequalities \eq{7x} and \eq{7y}, the condition {\it(iii)} follows (also in the case $x=u$).

{\it(iii)$\Rightarrow$(iv):} 
To deduce {\it(iv)} from {\it(iii)}, let $x_1,\dots,x_n\in I, t_1,\dots,t_n\geq 0$ with $t_1+\dots+t_n=1$ and $u:=t_1x_1+\dots+t_nx_n$.  Substituting $x$ by $x_i$ in the inequality of condition {\it(iii)}, then multiplying this inequality by $t_i$, finally adding up the inequalities so obtained side by side, we get
\Eq{*}{
\sum_{i=1}^{n} t_i\big(f(u)+(x_i-u)\varphi(u)\big)\leq&\sum_{i=1}^{n} t_i\big(f(x_i)+\Phi(|u-x_i|)\big).
}
Using that $\sum_{i=1}^{n} t_i(x_i-u)=0$, the above inequality simplifies to the inequality of condition {\it(iv)}.

{\it(iv)$\Rightarrow$(i):} To deduce the $\Phi$-convexity of $f$ from condition {\it(iv)}, let $x,y\in I$ and $t\in[0,1]$. Taking $n=2$, $x_1:=x$, $x_2:=y$, $t_1:=t$ and $t_2:=1-t$ in condition {\it(iv)}, it is immediate to see that the inequality reduces to the defining inequality of $\Phi$-convexity.
\end{proof}

Motivated by the condition {\it(iii)} of \thm{1}, we say that $\varphi:I\to\R$ is a \emph{$\Phi$-slope function for $f$} if it satisfies inequality \eq{ux} for all $x,u\in I$. Then \thm{1} states that the $\Phi$-convexity of a function $f$ is equivalent to the existence of a $\Phi$-slope function for $f$.

\Thm{1.2}{Let $\Phi\in\E(I)$ and $f:I\to\R$. Then the following conditions are equivalent to each other:
\begin{enumerate}[(i)]
 \item $f$ is $\Phi$-affine;
 \item For all $x,u,y\in I$ with $x<u<y$, 
 \Eq{*}{
\bigg|\dfrac{f(u)-f(x)}{u-x}-\dfrac{f(y)-f(u)}{y-u}\bigg|\leq\dfrac{\Phi(u-x)}{u-x} +\dfrac{\Phi(y-u)}{y-u};
 }
 \item Provided that $\Phi$ is increasing, there exists a function $\varphi: I\to\R$ such that, for all $x,u\in I$,
 \Eq{aslf}{
 |f(u)-f(x)-(u-x)\varphi(u)|\leq \Phi(|u-x|);
 }
 \item For all $n\in\N$, $x_1,\dots,x_n\in I,$ $t_1,\dots,t_n\geq0$ with $ t_1+\dots+t_n=1$,
 \Eq{*}{
 \Big|f(t_1x_1+\dots+t_nx_n)-\sum_{i=1}^{n}t_i\big(f(x_i)\big)\Big|\leq \sum_{i=1}^{n}t_i\Phi(|(t_1x_1+\dots+t_nx_n)-x_i|\big).
 }
\end{enumerate}
}

\begin{proof} Assertions (i), (ii) and (iv) are equivalent by \thm{1}.

Assume that (ii) holds and $\Phi$ is increasing.
First we are going to show the following inequality:
\Eq{minmax}{
  \max\bigg(&\sup_{x<u}\frac{f(u)-f(x)-\Phi(u-x)}{u-x},
     \sup_{u<y}\frac{f(y)-f(u)-\Phi(y-u)}{y-u}\bigg)\\
  &\leq \min\bigg(\inf_{x<u}\frac{f(u)-f(x)+\Phi(u-x)}{u-x},
     \inf_{u<y}\frac{f(y)-f(u)+\Phi(y-u)}{y-u}\bigg).
}
By condition (ii), for all $x,u,y\in I$ with $x<u<y$, we have the following two inequalities
\Eq{2I}{
  \dfrac{f(u)-f(x)}{u-x}-\dfrac{f(y)-f(u)}{y-u}
  &\leq\dfrac{\Phi(u-x)}{u-x} +\dfrac{\Phi(y-u)}{y-u}, \\
  -\dfrac{f(u)-f(x)}{u-x}+\dfrac{f(y)-f(u)}{y-u}
  &\leq\dfrac{\Phi(u-x)}{u-x} +\dfrac{\Phi(y-u)}{y-u}.
}
Then, one can easily see that the first and the second of the above inequalities 
imply
\Eq{*}{
  \sup_{x<u}\frac{f(u)-f(x)-\Phi(u-x)}{u-x}
  &\leq \inf_{u<y}\frac{f(y)-f(u)+\Phi(y-u)}{y-u}, \\
  \sup_{u<y}\frac{f(y)-f(u)-\Phi(y-u)}{y-u}
  &\leq \inf_{x<u}\frac{f(u)-f(x)+\Phi(u-x)}{u-x}
}
for all $u\in I$, respectively. In order to show that \eq{minmax} holds, it remains to prove that
\Eq{IRTP}{
  \sup_{x<u}\frac{f(u)-f(x)-\Phi(u-x)}{u-x}
  &\leq \inf_{x<u}\frac{f(u)-f(x)+\Phi(u-x)}{u-x},\\
  \sup_{u<y}\frac{f(y)-f(u)-\Phi(y-u)}{y-u}
  &\leq \inf_{u<y}\frac{f(y)-f(u)+\Phi(y-u)}{y-u}.
}
For the first inequality, we need to prove that, for all $x_1,x_2<u$, 
\Eq{A1}{
  \frac{f(u)-f(x_1)-\Phi(u-x_1)}{u-x_1}\leq \frac{f(u)-f(x_1)+\Phi(u-x_1)}{u-x_1}
}
If $x_1=x_2$, this follows from the nonnegativity of $\Phi$. We have to consider two cases according to the possibilites $x_1<x_2$ or $x_2<x_1$. In the first case, using \eq{2I} for the triplet $x_1<x_2<u$, then the monotonicity of $\Phi$, we get
\Eq{*}{
  \frac{f(u)-f(x_1)}{u-x_1}-\frac{f(u)-f(x_2)}{u-x_2}
  &=\frac{x_2-x_1}{u-x_1}\bigg(\frac{f(x_2)-f(x_1)}{x_2-x_1}-\frac{f(u)-f(x_2)}{u-x_2}\bigg)\\
  &\leq\frac{x_2-x_1}{u-x_1}\bigg(\frac{\Phi(x_2-x_1)}{x_2-x_1}+\frac{\Phi(u-x_2)}{u-x_2}\bigg)\\
  &\leq\frac{\Phi(x_2-x_1)}{u-x_1}+\frac{x_2-x_1}{u-x_1}\frac{\Phi(u-x_2)}{u-x_2}\leq\frac{\Phi(u-x_1)}{u-x_1}+\frac{\Phi(u-x_2)}{u-x_2},
}
which shows \eq{A1} in the case $x_1<x_2$. In the case $x_2<x_1$, the proof is analogous, and hence the first inequality in \eq{IRTP} has been verified. The second inequality in \eq{IRTP} can be shown in a similar manner. Thus, inequality \eq{minmax} has been proved, too. Define now the function $\varphi:I\to\R$ by
\Eq{*}{
  \varphi(u):=\max\bigg(&\sup_{x<u}\frac{f(u)-f(x)-\Phi(u-x)}{u-x},
     \sup_{u<y}\frac{f(y)-f(u)-\Phi(y-u)}{y-u}\bigg).
}
Then, for all $x<u<y$, we get
\Eq{AA}{
  f(u)-f(x)-(u-x)\varphi(u)\leq \Phi(u-x), \qquad
  f(y)-f(u)+(u-y)\varphi(u)\leq \Phi(y-u)
}
Based on inequality \eq{minmax}, we also have 
\Eq{*}{
\varphi(u)\leq \min\bigg(\inf_{x<u}\frac{f(u)-f(x)+\Phi(u-x)}{u-x},
     \inf_{u<y}\frac{f(y)-f(u)+\Phi(y-u)}{y-u}\bigg), 
}
which implies, for all $x<u<y$, that
\Eq{BB}{
  f(x)-f(u)+(u-x)\varphi(u)\leq \Phi(u-x), \qquad
  f(u)-f(y)-(u-y)\varphi(u)\leq \Phi(y-u).
}
The first inequalities in \eq{AA} and \eq{BB} yield, for all $x<u$, that 
\Eq{*}{
  |f(u)-f(x)-(u-x)\varphi(u)|\leq \Phi(u-x).
}
On the other hand, the second inequalities in \eq{AA} and \eq{BB} show, for all $u<y$, that 
\Eq{*}{
  |f(u)-f(y)-(u-y)\varphi(u)|\leq \Phi(y-u).
}
Combining these two inequalities, we can see that assertion (iii) must be valid.

Conversely, if assertion (iii) holds, then $\varphi$ is a $\Phi$-slope function for $f$, whence it follows that $f$ is $\Phi$-convex. Similarly, (iii) also implies that $(-\varphi)$ is a $\Phi$-slope function for $(-f)$, hence $(-f)$ is also $\Phi$-convex. These two $\Phi$-convexity properties yield that $f$ is $\Phi$-affine.
\end{proof}

Motivated by the condition {\it(iii)} of \thm{1.2}, we say that $\varphi:I\to\R$ is an \emph{absolute $\Phi$-slope function for $f$} if it satisfies inequality \eq{aslf} for all $x,u\in I$. Then \thm{1.2} states that a function $f$ is $\Phi$-affine (provided that $\Phi$ is increasing) if and only if there exists an absolute $\Phi$-slope function for $f$.

The next two propositions state that the $\Phi$-slope function of a $\Phi$-convex (resp.\ absolute $\Phi$-slope function of a $\Phi$-affine) function is approximately monotone (resp.\ H\"older) with respect to a transformation $\Phi^*$ of the error function $\Phi$.

\Prp{1}{Let $\Phi\in\E(I)$, $f:I\to\R$ be a $\Phi$-convex function and $\varphi:I\to\R$ be a $\Phi$-slope function for $f$. Then $\varphi$ is $\Phi^*$-monotone, where
\Eq{Phi*}{
  \Phi^*(t):=\begin{cases}
             2\dfrac{\Phi(t)}{t} &\mbox{if } 0<t<\ell(I),\\
             0 &\mbox{if } t=0.
             \end{cases}
}}

\begin{proof} Let $\varphi:I\to\R$ be a $\Phi$-slope function for $f$ and let $x,y\in I$ with $x<y$ be arbitrary. By obvious substitution into the inequality \eq{ux}, we get
\Eq{*}{
  f(x)+(y-x)\varphi(x)\leq f(y)+\Phi(y-x),\\
  f(y)+(x-y)\varphi(y)\leq f(x)+\Phi(y-x). 
}
Adding these inequalities side by side, the values $f(x)$ and $f(y)$ cancel out, hence we arrive at
\Eq{*}{
  (y-x)(\varphi(x)-\varphi(y))\leq 2\Phi(y-x).
}
Dividing by $y-x$, this inequality implies that $\varphi$ is $\Phi^*$-monotone, which was to be proved.
\end{proof}

\Prp{2}{Let $\Phi\in\E(I)$ be increasing, $f:I\to\R$ be a $\Phi$-affine function and $\varphi:I\to\R$ be a $\Phi$-slope function for $f$. Then $\varphi$ is $\Phi^*$-monotone, where $\Phi*$ is defined by \eq{Phi*}.}

\begin{proof} Let $\varphi:I\to\R$ be an absolute $\Phi$-slope function for $f$ and let $x,y\in I$ with $x\neq y$ be arbitrary. By obvious substitution into the inequality \eq{aslf}, we get
\Eq{*}{
  |f(x)-f(y)-(x-y)\varphi(x)|\leq\Phi(|x-y|),\\
  |f(y)-f(x)-(y-x)\varphi(y)|\leq\Phi(|y-x|). 
}
Therefore, by the triangle inequality,
\Eq{*}{
  |(y-x)(\varphi(x)-\varphi(y))|
  &\leq
  |f(x)-f(y)-(x-y)\varphi(x)|+|f(y)-f(x)-(y-x)\varphi(y)|\\
  &\leq \Phi(|y-x|)
}
Dividing by $|y-x|$, this inequality shows that $\varphi$ is $\Phi^*$-H\"older, which was to be proved.
\end{proof}

\section{Optimal error functions}

We say that an error function $\Phi\in\E(I)$ \emph{possesses the property $\Gamma$} if it satisfies the inequality
\Eq{G}{
  \Phi(x+y)\leq\Phi(x)+\frac{2x+y}{y}\Phi(y) \qquad(x\geq0,\,y>0,\,x+y<\ell(I)).
}
The subclass of error functions in $\E(I)$ with the property $\Gamma$ will be denoted by $\E^\Gamma(I)$. The subset of $\E^\Gamma(I)$ whose elements also satisfy $\Phi(0)=0$ will be denoted by $\E_0^\Gamma(I)$.

One can easily see that any subadditive error function satisfies property $\Gamma$, however, as we will see later, the reversed implication is not true.

The next result establishes some necessary and some sufficient conditions for the $\Gamma$ property.

\Thm{CGP}{Let $\Phi\in\E^\Gamma(I)$. Then $\sqrt{\Phi}$ and the map $t\mapsto t^{-1}\Phi(t)$ is subadditive on $\,]0,\ell(I)[\,$. If, in addition, $\varphi:[0,\ell(I)[\,\to\R_+$ is decreasing on $\,]0,\ell(I)[\,$, then $\varphi\cdot\Phi\in\E^\Gamma(I)$. In particular, if $\Psi\in\E(I)$ and $t\mapsto t^{-2}\Psi(t)$ is deccreasing on $\,]0,\ell(I)[\,$, then $\Psi\in\E^\Gamma(I)$.}

\begin{proof} To prove the subadditivity of $\sqrt{\Phi}$ on $\,]0,\ell(I)[\,$, let $x,y>0$ with $x+y<\ell(I)$. We have
\Eq{*}{
  \Phi(x+y)\leq\Phi(x)+\frac{2x+y}{y}\Phi(y)
  =\Phi(x)+\Phi(y)+2\frac{x}{y}\Phi(y).
}
and, interchanging the roles of $x$ and $y$, we also have
\Eq{*}{
  \Phi(x+y)\leq\Phi(x)+\Phi(y)+2\frac{y}{x}\Phi(x).
}
These two inequalities imply that
\Eq{300}{
  \Phi(x+y)
  \leq\Phi(x)+\Phi(y)+2\min\Big(\frac{x}{y}\Phi(y),\frac{y}{x}\Phi(x)\Big).
}
By using that the geometric mean of two numbers exceeds their minimum, we get
\Eq{*}{
  \Phi(x+y)\leq\Phi(x)+\Phi(y)+2\sqrt{\Phi(x)\Phi(y)}
  =\big(\sqrt{\Phi(x)}+\sqrt{\Phi(y)}\big)^2.
}
Taking square root of this inequality side by side, it follows that
\Eq{*}{
  \sqrt{\Phi(x+y)}\leq\sqrt{\Phi(x)}+\sqrt{\Phi(y)},
}
which completes the proof of the subadditivity of $\sqrt{\Phi}$.

To show the subadditivity of $t\mapsto t^{-1}\Phi(t)$, we use that the minmimu of two numbers is smaller than their arithmetic mean. Thus, from \eq{300} we arrive at
\Eq{*}{
  \Phi(x+y)\leq\Phi(x)+\Phi(y)+\Big(\frac{x}{y}\Phi(y)+\frac{y}{x}\Phi(x)\Big)=(x+y)\frac{\Phi(x)}{x}+(x+y)\frac{\Phi(y)}{y}.
}
Dividing both sides by $x+y$, we can reach at our desired conclusion.

To prove the next assertion, denote $\Psi:=\varphi\cdot\Psi$ and let $x>0$, $y>0$ with $x+y<\ell(I)$. Using the $\Gamma$ property of $\Phi$ and then the decreasingness of $\varphi$, we get
\Eq{*}{
  \Psi(x+y)=\varphi(x+y)\Phi(x+y)
  &\leq\varphi(x+y)\Phi(x)+\frac{2x+y}{y}\varphi(x+y)\Phi(y)\\
  &\leq\varphi(x)\Phi(x)+\frac{2x+y}{y}\varphi(y)\Phi(y)
  =\Psi(x)+\frac{2x+y}{y}\Psi(y),
}
which establishes the $\Gamma$ property of $\Psi$.

Finally, let $\Psi\in\E(I)$ be arbitrary such that $t\mapsto t^{-2}\Psi(t)=:\varphi(t)$ is decreasing on $\,]0,\ell(I)[\,$. Since $\Phi(t):=t^2$ possesses the $\Gamma$ property, it follows that $Ł\Psi=\varphi\cdot\Phi$ belongs to $\E^\Gamma(I)$.
\end{proof}

\Cor{pow}{Let $p\in\R$ and define the function $\Phi_p:[0,\infty[\,\to\R_+$ by $\Phi_p(t):=t^p$ for $t>0$ and $\Phi_p(0):=0$. Then $\Phi_p\in\E^\Gamma(\R_+)$ if and only if $p\leq 2$.}

\begin{proof} Assume that $\Phi_p$ possesses the $\Gamma$ property. Then, according to the first part of \thm{CGP}, $t\mapsto t^{-1}\Phi_p(t)=t^{p-1}$ is subadditive on $\R_+$. Therefore, $p-1\leq 1$, which yields $p\leq 2$.  

On the other hand, if $p\leq 2$, then $t\mapsto t^{-2}\Phi_p(t)=t^{p-2}$ is decreasing, therefore the last part of \thm{CGP} yields that $\Phi_p\in\E^\Gamma(\R_+)$.
\end{proof}

\Prp{3.1}{Let $\Phi\in\E^\Gamma(I)$. Then, for all $n\in\N$ and for all $u_1\geq0$, $u_2,\dots,u_n>0$ with $u_1+\dots+u_n<\ell(I)$, the following inequality holds
\Eq{3.1}{
  \Phi(u_1+\dots+u_n)\leq \Phi(u_1)+\frac{2u_1+u_2}{u_2}\Phi(u_2)
   +\dots+\frac{2(u_1+\dots+u_{n-1})+u_n}{u_n}\Phi(u_n).
}}

\begin{proof} Assume that $\Phi$ possesses the $\Gamma$ property, i.e., \eq{G} is satisfied. 

If $n=1$, then the statement trivially holds with equality. If $n=2$, then \eq{G} with $x:=u_1$ and $y:=u_2$ is equivalent to \eq{3.1}. Assume that \eq{3.1} has been proved for some $n\geq 2$ and let $u_1\geq0$, $u_2,\dots,u_{n+1}>0$ with $u_1+\dots+u_{n+1}<\ell(I)$. Then, by the $\Gamma$ property with $x:=u_1+\dots+u_n$, $y:=u_{n+1}$ and then by inductive assumption, we get
\Eq{*}{
  \Phi(u_1&+\dots+u_{n+1}) \\ 
  & \leq \Phi(u_1+\dots+u_n)+\frac{2(u_1+\dots+u_n)+u_{n+1}}{u_{n+1}}\Phi(u_{n+1}) \\
  &\leq \Phi(u_1)
   +\dots+\frac{2(u_1+\dots+u_{n-1})+u_n}{u_n}\Phi(u_n) 
   +\frac{2(u_1+\dots+u_n)+u_{n+1}}{u_{n+1}}\Phi(u_{n+1}),
}
which is exactly the inequality to be proved for $n+1$ variables.
\end{proof}

\Prp{3.19}{The classes $\E^\Gamma(I)$ and $\E_0^\Gamma(I)$ are closed with respect to addition, multiplication by nonnegative scalars, pointwise maximum, pointwise chain infimum and the limsup operation.}

\begin{proof} The proof is similar to that of \prp{MC1} and \prp{MC2}, therefore it is left to the reader.
\end{proof}

In view of the above proposition, we can see that the supremum of all error functions with the $\Gamma$ property which are smaller than a given error function $\Phi\in\E(I)$ is the largest error function below $\Phi$ with the $\Gamma$ property, which will be denoted as $\Phi^\Gamma$ and called the \emph{$\Gamma$-envelope of $\Phi$}. In what follows, we describe a construction for $\Phi^\Gamma$ in the class $\E_0^\Gamma(I)$. For this aim, for any error function $\Phi\in\E_0(I)$, we define its \emph{$\gamma$-transform} $\Phi^\gamma:[0,\ell(I)[\,\to\R_+$ of $\Phi$ by
\Eq{*}{
\Phi^\gamma(0):=0, \qquad
\Phi^\gamma(u):=\inf\bigg\{\Phi(x)+\frac{2x+y}{y}\Phi(y)\,
\bigg|\,x\geq0,\,y>0\colon x+y=u\bigg\} \qquad(u>0).  
}
Obviously, taking $x=0$ in the above definition, it follows from $\Phi(0)=0$ that $\Phi^\gamma\leq\Phi$ on $[0,\ell(I)[$.

\Thm{3.20}{
Let $\Phi\in\E_0(I)$ and define the sequence $(\Phi_n)_{n=0}^\infty$ by
\Eq{Git}{
  \Phi_1:=\Phi,\qquad \Phi_{n+1}:=\Phi_{n}^\gamma \qquad(n\in\N).
}
Then $(\Phi_n)$ is a pointwise decreasing sequence on $[0,\ell(I)[$ whose pointwise limit function equals $\Phi^\Gamma$ on $[0,\ell(I)[$.}

\begin{proof} The inequality $\Phi_2\leq\Phi_1$ is the consequence of $\Phi^\gamma\leq \Phi$. Now, using the obvious monotonicity of the mapping $\Phi\mapsto\Phi^\gamma$, by induction on $n$, it follows that $\Phi_{n+1}\leq\Phi_n$ holds on $[0,\ell(I)[$, which proves the decreasingness of the sequence $(\Phi_n)$. This implies that the pointwise limit function of this sequence, denoted as $\Phi_0$, exists and is not greater than $\Phi$. 

By the construction, for all $n\in\N$, we have that $\Phi_n(0)=0$ and
\Eq{*}{
\Phi_{n+1}(u)=\Phi_n^\gamma(u)=\inf\bigg\{\Phi_n(x)+\frac{2x+y}{y}\Phi_n(y)\,
\bigg|\,x\geq0,\,y>0\colon x+y=u\bigg\} \qquad(0<u<\ell(I)).  
}
This implies that
\Eq{*}{
  \Phi_{n+1}(x+y)\leq \Phi_n(x)+\frac{2x+y}{y}\Phi_n(y)
  \qquad(x\geq0,\,y>0,\, x+y<\ell(I)).
}
Upon taking the limit as $n\to\infty$, it follows that $\Phi_0$ possesses the $\Gamma$ property. 

It remains to show that $\Phi_0$ is nonsmaller than any error function $\Psi\in\E(I)$ with $\Psi\leq\Phi$ which possesses the $\Gamma$ property. We prove, by induction on $n$, that $\Psi\leq\Phi_n$ on  $[0,\ell(I)[$. The inequality $\Psi\leq\Phi_1=\Phi$ holds by assumption. Assume that $\Psi\leq\Phi_n$ is valid for some $n\in\N$. Let $u\in \,]0,\ell(I)[$ and let $x\geq0$, $y>0$ such that $u=x+y$.
Then
\Eq{*}{
  \Psi(u)=\Psi(x+y)
  \leq \Psi(x)+\frac{2x+y}{y}\Psi(y) 
  \leq \Phi_n(x)+\frac{2x+y}{y}\Phi_n(y).
}
Upon taking the infimum for all $x\geq0$, $y>0$ with $u=x+y$, it follows that
\Eq{*}{
  \Psi(u)\leq\Phi_n^\gamma(u)=\Phi_{n+1}(u).
}
This shows that $\Psi\leq\Phi_{n+1}$ holds on $[0,\ell(I)[$, which was to be verified.

Finally, taking the limit $n\to\infty$, we can conclude that $\Psi\leq\Phi_0$ is valid on $[0,\ell(I)[$, which demonstrates that $\Phi_0$ is the largest error function below $\Phi$ possessing the $\Gamma$ property. Therefore, $\Phi_0=\Phi^\Gamma$ holds.
\end{proof}

The next result establishes a sufficient condition that ensures the $\Gamma$ envelope to be zero.

\Thm{lim}{Let $\Phi\in\E(I)$ and assume that
\Eq{lim}{
  \lim_{t\to0^+} t^{-2}\Phi(t)=0.
}
Then $\Phi^\Gamma(u)=0$ holds for $u\in\,]0,\ell(I)[\,$.}

\begin{proof} Let $u\in\,]0,\ell(I)[\,$ and let $n\in\N$. Then, with $u_1:=\cdots:=u_n:=\frac{u}{n}$, the inequality \eq{3.1} applied for $\Phi^\Gamma$ (instead of $\Phi$) and the inequality $\Phi^\Gamma\leq\Phi$ yield that
\Eq{*}{
  \Phi^\Gamma(u)
  \leq (1+3+\dots+(2n-1))\Phi^\Gamma\Big(\frac{u}{n}\Big)
  =n^2\Phi^\Gamma\Big(\frac{u}{n}\Big)
  \leq n^2\Phi\Big(\frac{u}{n}\Big)
  =u^2\Big(\frac{n}{u}\Big)^2\Phi\Big(\frac{u}{n}\Big).
}
Upon taking the limit $n\to\infty$ and using the relation \eq{lim}, we get that $\Phi^\Gamma(u)=0$, which proves the statement.
\end{proof}

\Cor{lim}{Let $p\in\R$ and define $\Phi_p$ as in \cor{pow}. Then
\Eq{*}{
\Phi_p^\Gamma=\begin{cases}
              \Phi_p &\mbox{ if } p\leq2,\\ 0 &\mbox{ if } p>2.
             \end{cases}
}}

\begin{proof} If $p\leq 2$, then, by \cor{pow}, $\Phi_p$ possesses the $\Gamma$ property, and hence, $\Phi_p^\Gamma=\Phi_p$. In the case $p>0$,
we have that
\Eq{*}{
  \lim_{t\to0^+} t^{-2}\Phi_p(t)=\lim_{t\to0^+} t^{p-2}=0.
}
Therefore, by \thm{lim}, it follows that $\Phi_p^\Gamma=0$.
\end{proof}

\Thm{3.3}{Let $\Phi\in\E_0(I)$. Then
\Eq{*}{
  \C_\Phi(I)=\C_{\Phi^\Gamma}(I) \qquad\mbox{and}\qquad
  \A_\Phi(I)=\A_{\Phi^\Gamma}(I).
}}

\begin{proof} Due to the inequality $\Phi^\Gamma\leq\Phi$, it follows that 
\Eq{*}{
  \C_\Phi(I)\supseteq\C_{\Phi^\Gamma}(I) \qquad\mbox{and}\qquad
  \A_\Phi(I)\supseteq\A_{\Phi^\Gamma}(I).
} 
To prove the reversed inclusions, we will first verify that $\C_\Phi(I)\subseteq\C_{\Phi^\gamma}(I)$.

Let $f$ be a $\Phi$-convex function and let $x,y\in I$ with $x<y$ and $t\in\,]0,1[\,$. Define $u:=tx+(1-t)y$ and then choose $v\in[x,u[\,$ and $w\in\,]u,y]$ arbitrarily. Applying the $\Phi$-convexity of $f$ three times, we obtain
\Eq{*}{
  f(u)
  &\leq \frac{w-u}{w-v}\big(f(v)+\Phi(u-v)\big)
       +\frac{u-v}{w-v}\big(f(w)+\Phi(w-u)\big)\\
  &\leq \frac{w-u}{w-v}\bigg(\frac{u-v}{u-x}\big(f(x)+\Phi(v-x)\big)
       +\frac{v-x}{u-x}\big(f(u)+\Phi(u-v)\big)+\Phi(u-v)\bigg) \\
  &\qquad +\frac{u-v}{w-v}\bigg(\frac{y-w}{y-u}\big(f(u)+\Phi(w-u)\big)
       +\frac{w-u}{y-u}\big(f(y)+\Phi(y-w)\big)+\Phi(w-u)\bigg)\\
}
After a simple calculation, this reduces to the inequality
\Eq{vw}{
  f(u)&\leq \frac{y-u}{y-x}\bigg(f(x)+\Phi(v-x)
     +\frac{u+v-2x}{u-v}\Phi(u-v)\bigg)\\
    &\qquad +\frac{u-x}{y-x}\bigg(f(y)+\Phi(y-w)+\frac{2y-u-w}{w-u}\Phi(w-u)\bigg).
}
Observe, that by the definition of $\Phi^\gamma$, we have that
\Eq{*}{
  \inf_{v\in[x,u[}
  \bigg(\Phi(v-x)+\frac{u+v-2x}{u-v}\Phi(u-v)\bigg)=\Phi^\gamma(u-x)
}
and
\Eq{*}{
  \inf_{w\in\,]u,y]}
  \bigg(\Phi(y-w)+\frac{2y-u-w}{w-u}\Phi(w-u)\bigg)=\Phi^\gamma(y-u).
}
Therefore, upon taking the infimum with respect to $v\in[x,u[\,$ and $w\in\,]u,y]$ in \eq{vw}, it follows that
\Eq{*}{
  f(u)\leq \frac{y-u}{y-x}\big(f(x)+\Phi^\gamma(u-x)\big)
    +\frac{u-x}{y-x}\big(f(y)+\Phi^\gamma(y-u)\big),
}
which proves that $f$ is $\Phi^\gamma$-convex. 

Now define the sequence of error function $\Phi_n$ by the iteration \eq{Git}. 
By the assumption, $f$ is $\Phi_1=\Phi$-convex. On the other hand, using what we have proved above, it follows that if $f$ is $\Phi_n$-convex, then it is also $\Phi_{n+1}=\Phi_n^\gamma$-convex. Therefore, $f$ is $\Phi_n$-convex for all $n\in\N$. Upon taking the limit $n\to\infty$ and using \thm{3.20}, we obtain that $f$ is $\Phi^\Gamma$-convex.

For the second part of the assertion, observe that if $f$ is $\Phi$-affine, then by \prp{H0}, $f$ and $(-f)$ are $\Phi$-convex. Thus, by the first part, we get that $f$ and $(-f)$ are $\Phi^\Gamma$-convex, which according to \prp{H0} again implies that $f$ is $\Phi^\Gamma$-affine. Therefore, the reversed inclusion $\A_\Phi(I)\subseteq\A_{\Phi^\Gamma}(I)$ also holds.
\end{proof}

\Cor{lim+}{Let $p>2$ and define $\Phi_p$ as in \cor{pow}. Then
\Eq{*}{
  \C_{\Phi_p}(I)=\C_0(I) \qquad\mbox{and}\qquad
  \A_{\Phi_p}(I)=\A_0(I).
}}

\begin{proof} 
The statement directly follows from \thm{3.3} and \cor{lim}.
\end{proof}

\section{The $\Phi$-convex envelop}

For $\Phi\in\E_0(I)$ and $f:I\to\R$, define $\underline{C}_\Phi(f):I\to[-\infty,\infty)$ by
\Eq{***}{
  \underline{C}_\Phi(f)(u)
  :=\inf_{\mbox{\tiny$\begin{array}{c}x,y\in I,\,t\in[0,1]\\u=tx+(1-t)y\end{array}$}}t\big(f(x)+\Phi((1-t)|y-x|)\big)+(1-t)\big(f(y)+\Phi(t|y-x|)\big).
}

\Thm{CPhi}{Let $\Phi\in\E_0(I)$. Then the operator $\underline{C}_\Phi$ is monotone and concave in the pointwise sense, i.e., for $f,g:I\to\R$ and $\lambda\in[0,1]$, we have
\Eq{*}{
  \lambda \underline{C}_\Phi(f)+(1-\lambda)\underline{C}_\Phi(g)
  \leq \underline{C}_\Phi(\lambda f+(1-\lambda)g).
}
Furthermore, for all $f:I\to\R$, the inequality $\underline{C}_\Phi(f)\leq f$ holds and here equality is valid if and only if $f$ is $\Phi$-convex. In other words, the fixed points of $\underline{C}_\Phi$ are exactly the $\Phi$-convex functions.}

\begin{proof} The monotonicity of $\underline{C}_\Phi$ is obvious. To see its subadditivity, observe that, by the definition, for all $u\in I$ and $x,y\in I,\,t\in[0,1]$ with $u=tx+(1-t)y$,
we have
\Eq{*}{
   \lambda &\underline{C}_\Phi(f)(u)+(1-\lambda)\underline{C}_\Phi(g)(u)\\
   &\leq \lambda\big(t\big(f(x)+\Phi((1-t)|y-x|)\big)+(1-t)\big(f(y)+\Phi(t|y-x|)\big)\big) \\
   &\quad +(1-\lambda)\big(t\big(g(x)+\Phi((1-t)|y-x|)\big)+(1-t)\big(g(y)+\Phi(t|y-x|)\big)\big)\\
   &=t\big((\lambda f+(1-\lambda)g)(x)+\Phi((1-t)|y-x|)\big)+(1-t)\big((\lambda f+(1-\lambda)g)(y)+\Phi(t|y-x|)\big).
}
Upon taking the infimum for $x,y\in I,\,t\in[0,1]$ with $u=tx+(1-t)y$,
it follows that
\Eq{*}{
  \lambda \underline{C}_\Phi(f)(u)+(1-\lambda)\underline{C}_\Phi(g)(u)
  \leq \underline{C}_\Phi(\lambda f+(1-\lambda)g)(u).
}
Taking $x:=y:=u$ in the definition of $\underline{C}_\Phi(f)(u)$ and using $\Phi(0)=0$, we can see that $\underline{C}_\Phi(f)(u)\leq f(u)$ holds for all $u\in I$. 

On the other hand, if $f=\underline{C}_\Phi(f)$ holds for some $f:I\to\R$, then, for all $x,y\in I$ and $t\in[0,1]$ with $u:=tx+(1-t)y$, we get
\Eq{*}{
  f(tx+(1-t)y)&=f(u)=\underline{C}_\Phi(f)(u)\\
  &\leq t\big(f(x)+\Phi((1-t)|y-x|)\big)+(1-t)\big(f(y)+\Phi(t|y-x|)\big),
}
which shows that $f$ is $\Phi$-convex.

Conversely, if $f$ is $\Phi$-convex, then, for all $x,y\in I$ and $t\in[0,1]$ with $u=tx+(1-t)y$,
\Eq{*}{
  f(u)=f(tx+(1-t)y)\leq t\big(f(x)+\Phi((1-t)|y-x|)\big)+(1-t)\big(f(y)+\Phi(t|y-x|)\big).
}
Now, taking the infimum for $x,y\in I,\,t\in[0,1]$ with $u=tx+(1-t)y$,
it follows that $f(u)\leq \underline{C}_\Phi(f)(u)$. This implies that, in fact $f(u)=\underline{C}_\Phi(f)(u)$ an completes the proof of the equality $f=\underline{C}_\Phi(f)$.
\end{proof}

In what follows, we construct the $\Phi$-convex envelope of any function which admits a $\Phi$-convex minorant in terms of the operator $\underline{C}_\Phi$. The following auxiliary result will be useful.

\Lem{CCC}{Let $\Phi\in\E_0(I)$ be increasing and subadditive. Then $\Psi(u):=-\Phi(|u|)$ is $\Phi$-convex on $J:=I-I$. In particular, for $p\in[0,1]$, $\Psi_p(u):=-\Phi_p(|u|)$ is $\Phi_p$-convex on $\R$.}

\begin{proof} Let $x,y\in J$ and $t\in[0,1]$. Then
\Eq{*}{
  -\Psi(x)=\Phi(|x|)&\leq\Phi(|tx+(1-t)y|+|(1-t)(x-y)|)\\
  &\leq\Phi(|tx+(1-t)y|)+\Phi(|(1-t)(x-y)|)\\
  &=-\Psi(tx+(1-t)y)+\Phi(|(1-t)(x-y)|)
}
and
\Eq{*}{
  -\Psi(y)=\Phi(|y|)&\leq\Phi(|tx+(1-t)y|+|t(y-x)|)\\
  &\leq\Phi(|tx+(1-t)y|)+\Phi(|t(x-y)|)\\
  &=-\Psi(tx+(1-t)y)+\Phi(|t(x-y)|).
}
Multiplying the first inequality by $t$ and the second one by $(1-t)$ and then adding up the inequalities so obtained side by side, we obtain that
\Eq{*}{
  -t\Psi(x)-(1-t)\Psi(y)
  \leq -\Psi(tx+(1-t)y)+t\Phi(|(1-t)(x-y)|)+(1-t)\Phi(|t(x-y)|),
}
which proves that $\Psi$ is $\Phi$-convex.

The last assertion of the lemma follows from the fact the $\Phi_p$ is nondecreasing and subadditive provided that $p\in[0,1]$. 
\end{proof}

\Thm{mnv}{Let $\Phi\in\E_0(I)$ such that $\Psi(u):=-\Phi(|u|)$ is $\Phi$-convex on $J:=I-I$. Let $f:I\to \R$ be a function which admits a $\Phi$-convex minorant. Then the function $\underline{C}_\Phi(f)$ is the largest $\Phi$-convex function which is smaller than or equal to $f$.}

\begin{proof}
The inequality $\underline{C}_\Phi(f)\leq f$ follows from the previous theorem.

Now suppose $g$ is a $\Phi$-convex function such that $g\leq f$ holds (by the assumption, there is at least one such function $g$). Then, according to \thm{CPhi}, $g=\underline{C}_\Phi(g)\leq \underline{C}_\Phi(f)$, which proves that $g\leq \underline{C}_\Phi(f)$.

To complete the proof, it will be sufficient to show that $g:=\underline{C}_\Phi(f)$ is $\Phi$-convex. Let $u\in I$ be fixed arbitrarily. Then, for all $x<u<y$, with $t:=\frac{y-u}{y-x}$, we have $u=tx+(1-t)y$. Thus, the definition of $g(u)=\underline{C}_\Phi(f)(u)$ yields
\Eq{*}{
  g(u)\leq\frac{y-u}{y-x}\big(f(x)+\Phi(u-x)\big)
    +\frac{u-x}{y-x}\big(f(y)+\Phi(y-u)\big).
}
This inequality implies
\Eq{*}{
  \varphi(u):=\sup_{x<u}\frac{g(u)-f(x)-\Phi(u-x)}{u-x}
    \leq\inf_{u<y}\frac{f(y)+\Phi(y-u)-g(u)}{y-u}.
}
Therefore, for every $z\in I$, we have
\Eq{uz}{
   g(u)+\varphi(u)(z-u)\leq f(z)+\Phi(|z-u|).
}
Let $v\in I$ and $\varepsilon>0$ be arbitrary. Then, by the definition of $g(v)=\underline{C}_\Phi(f)(v)$, there exist $x,y\in I$, $t\in[0,1]$ with $v=tx+(1-t)y$ such that
\Eq{*}{
  t\big(f(x)+\Phi((1-t)|y-x|)\big)+(1-t)\big(f(y)+\Phi(t|y-x|)\big)
  <g(v)+\varepsilon.
}
Now, applying the inequality \eq{uz} for $z:=x$ and $z:=y$, and in the last step using the $\Phi$-convexity of $\Psi$, it follows that
\Eq{*}{
   g(u)&+\varphi(u)(v-u)\\
   &=t\big(g(u)+\varphi(u)(x-u)\big)+(1-t)\big(g(u)+\varphi(u)(y-u)\big)\\
   &\leq t\big(f(x)+\Phi(|x-u|)\big)+(1-t)\big(f(y)+\Phi(|y-u|)\big)\\
   &<g(v)+\varepsilon+t\big(\Phi(|x-u|)-\Phi((1-t)|y-x|)\big) +(1-t)\big(\Phi(|y-u|)-\Phi(t|y-x|)\big)\\
   &=g(v)+\varepsilon+t\big(-\Psi(x-u)-\Phi((1-t)|y-x|)\big) +(1-t)\big(-\Psi(y-u)-\Phi(t|y-x|)\big)\\
   &\leq g(v)+\varepsilon-\Psi((tx+(1-t)y)-u)
   = g(v)+\varepsilon+\Phi(|v-u|).
}
Upon taking the limit $\varepsilon\to0^+$, we obtain that $g$
satisfies the inequality
\Eq{*}{
   g(u)+\varphi(u)(v-u)\leq g(v)+\Phi(|v-u|) 
   \qquad(u,v\in I).
}
Therefore, according to the characterization theorem of $\Phi$-convexity, this implies that $g$ is $\Phi$-convex.
\end{proof}

Combining the assertions of \lem{CCC} and \thm{mnv}, we immediately obtain the following consequence. 

\Cor{mnv}{Let $\Phi\in\E_0(I)$ be nondecreasing and subadditive. Let $f:I\to \R$ be a function which admits a $\Phi$-convex minorant. Then the function $\underline{C}_\Phi(f)$ is the largest $\Phi$-convex function which is smaller than or equal to $f$.}

The next result is an extension of the ...

\Thm{mnz}{Let $\Phi\in\E_0(I)$ and let $f,g:I \to \R$ such that $g \leq h \leq f$ for some $\Phi$-convex function $h: I\to \R$.
Then, for all $x,y \in I$ and $t \in [0,1]$, the functional inequality 
\Eq{550}
{
g(tx+(1-t)y)\leq tf(x)+(1-t)f(y)+t(\Phi(1-t)|y-x|)+(1-t)\Phi(t|y-x|)
}
holds. Conversely, if $\Psi(u):=-\Phi(|u|)$ is $\Phi$-convex on $J:=I-I$ and \eq{550} is valid for all $x,y \in I$ and $t \in [0,1]$, then 
there exists a $\Phi$-convex function $h: I \to \R$ such that $g \leq h \leq f$ holds on $I$.
}

\begin{proof}
Let $f,g: I\to \R$ such that $g\leq h\leq f$ holds for some $\Phi$-convex function $h$. Then, for all $x,y \in I$ and $t \in [0,1]$, we have
\Eq{*}
{
g(tx+(1-t)y)&\leq h(tx+(1-t)y)\\
&\leq th(x)+(1-t)h(y)+t(\Phi(1-t)|y-x|)+(1-t)\Phi(t|y-x|)\\
&\leq tf(x)+(1-t)f(y)+t(\Phi(1-t)|y-x|)+(1-t)\Phi(t|y-x|).
}
This shows the validity of \eq{550} on the indicated domain.

Assume that $f,g:I \to \R$ satisfy \eq{550}. Then, this inequality yields that $g\leq \underline{C}_\Phi(f)$. In view of \thm{mnv} and the $\Phi$-convexity of $\Psi$, it follows that $\underline{C}_\Phi(f)$ is $\Phi$-convex. Therefore, the statement is fulfilled with $h:=\underline{C}_\Phi(f)$.
\end{proof}

Combining the assertions of \lem{CCC} and \thm{mnz}, we immediately obtain the following consequence. 

\Cor{mnz}{Let $\Phi\in\E_0(I)$ be nondecreasing, subadditive and let $f,g:I \to \R$. Then the inequalities $g \leq h \leq f$ hold for some $\Phi$-convex function $h: I\to \R$ if and only if, for all $x,y \in I$ and $t \in [0,1]$, the functional inequality \eq{550} is satisfied.}

\Thm{mnv+}{Let $\Phi\in\E_0(I)$. Let $f:I\to \R$ be a function which admits a $\Phi$-convex minorant and define the sequence $f_n:I\to\R$ by
\Eq{*}{
   f_1:=f, \qquad f_{n+1}:=\underline{C}_\Phi(f_n) \qquad(n\in\N).
}
Then the sequence $(f_n)$ is pointwise decreasing and its limit function is the largest $\Phi$-convex function which is smaller than or equal to $f$.}

\begin{proof}
The pointwise monotonicity of the sequence follows from the property $\underline{C}_\Phi(f)\leq f$ established in \thm{CPhi}. Let $f_0$ denote  the pointwise limit function of the sequence $(f_n)$. Now suppose $g$ is an arbitrary $\Phi$-convex function such that $g\leq f$ Then $g\leq\underline{C}_\Phi(g)\leq\underline{C}_\Phi(f)=\underline{C}_\Phi(f_1)=f_2$. Using the same iterative argument, we can show that $g\leq f_n$ for all $n\in \N.$ Upon taking the limit, we get that $g\leq f_{0}$ holds. 

Finally, we show the $\Phi$-convexity of $f_0$. Let $x,y\in I$ and $t\in[0,1]$ be arbitrary. From definition of $f_{n+1}$ it follows that
\Eq{*}{
  f_{n+1}(tx+(1-t)y)
  \leq t\big(f_{n}(x)+\Phi((1-t)|y-x|)\big)+(1-t)\big(f_{n}(y)+\Phi(t|y-x|)\big).
}
Upon taking the limit as $n\to\infty$, we get that
\Eq{*}{
  f_{0}(tx+(1-t)y)
  \leq t\big(f_{0}(x)+\Phi((1-t)|y-x|)\big)+(1-t)\big(f_{0}(y)+\Phi(t|y-x|)\big),
}
which shows that $f_0$ is $\Phi$-convex, indeed.
\end{proof}


\end{document}